\documentclass[12pt]{amsart}
\usepackage{amsfonts,amsthm,amscd}

\newtheorem{theorem}{Theorem}[section]
\newtheorem{corollary}[theorem]{Corollary}
\newtheorem{lemma}[theorem]{Lemma}
\newtheorem{proposition}[theorem]{Proposition}
\newtheorem{remark}[theorem]{Remark}
\newtheorem{mtheorem}[theorem]{Main Theorem}
\theoremstyle{definition}
\newtheorem{defn}{Definition}
\theoremstyle{remark}
\newtheorem{exmp}{Example}
\newtheorem*{acknow}{Acknowledgments}
\newtheorem{case}{Case}
\newcommand{\nc}{\newcommand}
\newcommand{\be}{\begin{equation}}
\newcommand{\ee}{\end{equation}}
\newcommand{\bea}{\begin{eqnarray}}
\newcommand{\eea}{\end{eqnarray}}
\newcommand{\no}{\nonumber}

\newcommand{\bc}{\begin{center}}
\newcommand{\ec}{\end{center}}

\nc{\bmth}{\begin{mtheorem}}
\nc{\emth}{\end{mtheorem}}
\nc{\bth}{\begin{theorem}}
\nc{\eth}{\end{theorem}}
\nc{\bpr}{\begin{proposition}}
\nc{\epr}{\end{proposition}}
\nc{\ble}{\begin{lemma}}
\nc{\ele}{\end{lemma}}
\nc{\bco}{\begin{corollary}}
\nc{\eco}{\end{corollary}}
\nc{\bre}{\begin{remark}}
\nc{\ere}{\end{remark}}
\newcommand{\la}{\lambda}

\nc{\f}{\frac}
\nc{\rw}{\rightarrow}
\nc{\Rw}{\Rightarrow}
\newcommand{\st}[1]{\hat#1}
\nc{\nun}{\nabla ^{1} }
\nc{\ndo}{\nabla ^{2} }
\nc{\nt}{\stackrel{\sim}{\nabla}}
\nc{\De}{\Delta}

\nc{\na}{\nabla}
\nc{\al}{\alpha}
\nc{\rai}{R(e_\al ,e_i ,e_\al ,e_i)}
\nc{\bet}{\beta}
\nc{\rab}{R_{\al \bet \al \be}}
\nc{\ga}{\gamma}
\nc{\de}{\delta}
\nc{\ep}{\varepsilon}
\nc{\ei}{\ep _i}
\nc{\ea}{\ep _{\al}}
\nc{\ej}{\ep _j}
\nc{\eb}{\ep _{\bet}}
\nc{\br}{\bgroup\em}
\nc{\er}{\egroup}
\nc{\rs}{Riemannian submersion}
\nc{\srs}{semi-Riemannian submersion}
\nc{\co}{compact and orientable total space}
\nc{\com}{compact and orientable manifold}
\nc{\om}{\omega}
\nc{\ta}{\tau ^{HV}}

\newcommand{\pis}{$\pi :(M,g) \to (B,g')$}
\nc{\hauu}{h^\al _{11}}
\nc{\haij}{h^\al _{ij}}
\nc{\aiab}{A^i _{\al \bet}}
\nc{\dvg}{\mathrm{dv}_g}
\nc{\ve}{$\mathcal V$}
\nc{\ho}{$\mathcal H$}
\nc{\tm}{\Gamma (TM)}
\nc{\tb}{\Gamma (TB)}
\nc{\lc}{Levi-Civita connection}
\nc{\dec}{\st{\de}}
\nc{\deu}{\stackrel{\vee}{\de}}
\begin{document}
\title[Semi-Riemannian submersions from pseudo-hyperbolic
   spaces]{Semi-Riemannian submersions from real and complex
 pseudo-hyperbolic spaces}
%\author{Gabriel B\u adi\c toiu \and Stere Ianu\c s}
\author{Gabriel B\u adi\c toiu}
\author{Stere Ianu\c s}
\address{                            
    Institute of Mathematics             
    of the Romanian Academy,              
    P.O. Box 1-764,                      
    Bucharest 70700,                      
    Romania}
\email{gbadit@stoilow.imar.ro}
\address{University of Bucharest,
    Faculty of Mathematics, C.P. 10-119, Post. Of. 10,
    Bucharest 72200,
    Romania}
\email{ianus@geometry.math.unibuc.ro}
\date{}
    

\begin{abstract}
        We classify the semi-Riemannian submersions
        from a pseudo-hyper\-bolic space onto a Riemannian
        manifold under the assumption that the fibres 
        are connected and totally geodesic. Also we obtain
        the classification of the semi-Riemannian
        submersions from a complex pseudo-hyperbolic space onto a Riemannian
        manifold under the assumption that the fibres are complex, connected
        and totally geodesic submanifolds.  
\end{abstract}
\maketitle
\section*{Introduction}

The theory of Riemannian submersions was initiated 
by O'Neill \cite{one} and Gray \cite{gra}.
Presently, there is an extensive literature on the \rs s 
with different conditions imposed on the total
space and on the fibres.
A systematic exposition could be found in Besse's book \cite{bes} .
Semi-Riemannian submersions were introduced by
O'Neill in his book \cite{onei}.\\
The class of harmonic Riemannian submersions, 
and in particular of those with totally geodesic fibres,
is contained in the class of horizontally homothetic harmonic morphisms. For important
results concerning the geometry of harmonic morphisms we refer to \cite{bai}.
%In \cite{woo}, 
Wood constructs examples of harmonic morphisms from Riemannian
submersions with totally geodesic fibres by horizontally conformal deformation of the
metric. 
Recently, Fuglede studied harmonic morphisms between semi-Riemannian manifolds
(see \cite{fug}). In this paper we solve the classification problem of the
semi-Riemannian submersions with totally geodesic fibres from real and
complex pseudo-hyperbolic spaces.\\
R. Escobales \cite{esc}, \cite{esco} and A. Ranjan
\cite{ran} classified
Riemannian submersions with totally geodesic fibres
from a sphere $S^n$ and from a complex
projective space $\mathbb CP^n$. 
M.A. Magid \cite{mag} classified the
semi-Riemannian submersions with totally geodesic fibres from
an anti-de Sitter space onto a Riemannian manifold.
In section \S 2 we classify the semi-Riemannian submersions 
with totally geodesic fibres from a pseudo-hyperbolic space onto a 
Riemannian  manifold. Also we obtain the classification of
 the semi-Riemannian submersions 
with connected, complex, totally geodesic fibres
from a complex pseudo-hyperbolic space onto a Riemannian manifold.
%%%%%%%%%%%%%%%%%%%%%%%%%%%%%%%%%%%%%%%%%%%%%%%%%%%%%%%%%%%%%%%%%%%%%%%%%%%
\section{Preliminaries and examples}
\begin{defn}
Let $(M,g)$ be an $m$-dimensional connected 
semi-Riemannian manifold of index $s$ ($0 \le s\le m$),
let $(B,g')$  be an $n$-dimensional connected 
semi-Riemannian manifold of index $s'\le s$, ($0 \le s'\le n$).
A {\it semi-Riemannian submersion} (see \cite{onei}) is 
a smooth map $\pi :M\to B$
which is onto and satisfies the following three axioms:
\begin{itemize}
\item[(a)] $\pi_*|_p$ is onto for all $p\in M$;
\item[(b)] the fibres $\pi^{-1}(b)$, $b\in B$ are 
    semi-Riemannian submanifolds of $M$;
\item[(c)] $\pi_*$ preserves scalar products of 
    vectors normal to fibres.
\end{itemize}
\end{defn}
  We shall always assume that the dimension of the fibres
  $\dim M-\dim B$ is positive and
  the fibres are connected.

The tangent vectors to fibres are called 
vertical and those normal to fibres
are called horizontal.
We denote by \ve\ the vertical distribution and by \ho\
the horizontal distribution.\\
B. O'Neill \cite{one} has characterized the geometry 
of a Riemannian submersion in terms of
the tensor fields $T$, $A$ defined by
  \begin{eqnarray*}  A_E F&=&h \na_{hE}{vF} + 
                v \na_{hE}{hF}\\
   T_E F&=&h \na _{vE}{vF} + v \na_{vE}{hF}
  \end{eqnarray*}
for every $E$, $F$ tangent vector fields to $M$.
Here  $\na $ is the \lc\ of $g$, the symbols $v$ and $h$ are the orthogonal 
projections on \ve\ and \ho ,
respectively. 
The letters $U$, $V$ will always denote 
vertical vector fields,
 $X$, $Y$, $Z$ horizontal vector fields.
Notice that $T_UV$ is the second fundamental form
of each fibre and $A_XY$ is a natural {\it obstruction} to integrability of
horizontal distribution $\mathcal H$.
The tensor $A$ is called O'Neill's integrability tensor.
For basic properties of \rs s and examples see
\cite{bes}, \cite{gra}, \cite{one}.
A vector field $X$ on $M$ is said to be {\it basic} if 
$X$ is horizontal and $\pi -$related to a vector 
field $X'$ on $B$.
Notice that every vector field $X'$
on $B$ has a unique horizontal lifting $X$ to $M$ 
and $X$ is basic. 
The following lemma is well known (see \cite{one}).
\ble\label{lemdoi}
We suppose X and Y are basic vector fields on M
 which are $\pi$-related to $X'$ and
 $Y'$. Then
\begin{itemize} 
\item[a)] $h[X,Y]$ is basic and
 $\pi$-related to $[X',Y']$;
\item[b)] $h\na_X Y$ is basic and $\pi$-related to 
$\na '_{X'}Y'$\ , where $\na '$ is the \lc\ on B;
\end{itemize}
\ele

The O'Neill's integrability tensor $A$ has the following properties 
(see \cite{one} or \cite{bes}).
\ble\label{l:1}
Let $X$, $Y$ be horizontal vector fields and $E$, $F$ be vector fields on $M$.
Then each of the following holds:
\begin{itemize}
\item[(a)] $A_XY=-A_YX$;
\item[(b)] $A_{hE}F=A_{E}F$;
\item[(c)] $A_E$ maps the horizontal subspace into the vertical one
           and the vertical subspace into the horizontal one;
\item[(d)] $g(A_XE,F)=-g(E,A_XF)$;
\item[(e)] If moreover $X$ is basic then $A_XV=h\na_VX$ for every vertical vector
           field $V$;
\item[(f)] $g((\na_YA)_XE,F)=g(E,(\na_YA)_XF)$.
\end{itemize}
\ele

Let $\hat g$ be the induced metric on fibre $\pi^{-1}(\pi(p))$, $p\in M$.
We denote by $R$, $R'$, $\hat R$ the Riemann tensors
of the metrics $g$, $g'$, $\hat g$ respectively.

The following equations, usually called {\it O'Neill's equations}, 
characterize the geometry of a semi-Riemannian submersion 
(see \cite{one}, \cite{gra}, \cite{bes}).

\bpr
For every vertical vector fields $U$, $V$, $W$, $W'$ 
and for every horizontal vector fields $X$, $Y$, $Z$, $Z'$,
we have the following formulae:
\begin{itemize}
\item[$i)$] $R(U,V,W,W')=\hat R(U,V,W,W')-g(T_UW,T_VW')+g(T_VW,T_UW')$,
\item[$ii)$] $R(U,V,W,X)=g((\na _VT)_UW,X)-g((\na _UT)_VW,X),$
\item[$iii)$] $R(X,U,Y,V)=g((\na _XT)_UV,Y)-g(T_UX,T_VY)+g((\na
    _UA)_XY,V)+g(A_XU,A_YV)$,
\item[$iv)$] $R(U,V,X,Y)=g((\na _UA)_XY,V)-g((\na
    _VA)_XY,U)+g(A_XU,A_YV)-g(A_XV,A_YU)-g(T_UX,T_VY)+g(T_VX,T_UY)$,
\item[$v)$]  $R(X,Y,Z,U)=g((\na_ZA)_XY,U)+g(A_XY,T_UZ)-g(A_YZ,T_UX)-g(A_ZX,T_UY)$,
\item[$vi)$] $R(X,Y,Z,Z')=R'(\pi_*X,\pi_*Y,\pi_*Z,\pi_*Z')
	-2g(A_XY,A_ZZ')+g(A_YZ,A_XZ')-g(A_XZ,A_YZ').$
\end{itemize}
\epr

Using O'Neill's equations, we get the following lemma.
\ble\label{lemtrei}
If \pis\ is a \srs\ with totally geodesic fibres then:
\bea
a)\ R(U,V,U,V)&=&\st{R}(U,V,U,V);\no\\
b)\ R(X,U,X,U)&=&g(A_XU,A_XU);\no\\
c)\ R(X,Y,X,Y)&=&R'(\pi_*X,\pi_*Y,\pi_*X,\pi_*Y)-3g(A_XY,A_XY).\no
\eea
\ele

We recall the definitions of real and complex pseudo-hyperbolic spaces (see \cite{onei}
and \cite{bar}).
\begin{defn}
   Let $<\cdot,\cdot>$ be the symmetric bilinear form on $\mathbb R^{m+1}$
   given by
        $$<x,y>=-\sum\limits_{i=0}^sx_iy_i+\sum\limits_{i=s+1}^mx_iy_i$$
  for $x=(x_0,\cdot\cdot\cdot,x_m),y=(y_0,\cdot\cdot\cdot,y_m)\in\mathbb R^{m+1}$.
  For $s>0$ let\\
  $H^m_s=\{x\in\mathbb R^{m+1}\ |\ <x,x>=-1\}$ be the semi-Riemannian submanifold
  of 
    $$\mathbb R_{s+1}^{m+1}=(\mathbb R^{m+1}, 
  ds^2=-dx^0\otimes dx^0-\cdot\cdot\cdot -dx^s\otimes dx^s+
        dx^{s+1}\otimes dx^{s+1}+\cdot\cdot\cdot +dx^m\otimes dx^m).$$ 
  $H^m_s$ is called the $m$-dimensional ({\it real}) {\it pseudo-hyperbolic space} of index $s$. 
  We notice that $H^m_s$ has constant sectional curvature $-1$ and the curvature tensor is given by
  $R(X,Y,X,Y)=-g(X,X)g(Y,Y)+g(X,Y)^2$. $H^m_s$ can be written as homogeneous space, namely we have
  $H^m_s=SO(s+1,m-s)/SO(s,m-s)$, $H^{2m+1}_{2s+1}=SU(s+1,m-s)/SU(s,m-s)$,
  $H^{4m+3}_{4s+3}=Sp(s+1,m-s)/Sp(s,m-s)$ (see \cite{wol}).
\end{defn}
\begin{defn}\label{d:3}
   Let $(\cdot,\cdot)$ be the hermitian scalar product on $\mathbb C^{m+1}$ given by
     $$(z,w)=-\sum\limits_{i=0}^sz_i\bar{w_i}+\sum\limits_{i=s+1}^mz_i\bar{w_i}$$
  for $z=(z_0,\cdot\cdot\cdot,z_m), w=(w_0,\cdot\cdot\cdot,w_m)\in\mathbb C^{m+1}$.
  Let $M$ be the real hypersurface of $\mathbb C^{m+1}$
  given by $M=\{z\in\mathbb C^{m+1}\ |\ (z,z)=-1\}$ and
  endowed with the induced metric of
  $$(\mathbb C^{m+1},
  ds^2=-dz^0\otimes d\bar{z}^0-\cdot\cdot\cdot -dz^s\otimes d\bar{z}^s+
        dz^{s+1}\otimes d\bar{z}^{s+1}+\cdot\cdot\cdot +dz^m\otimes d\bar{z}^m).$$ 
  The natural action of $S^1=\{e^{i\theta}\ |\ \theta\in\mathbb R\}$ on $\mathbb C^{m+1}$
  induces an action on $M$. Let $\mathbb CH^m_s=M/S^1$ endowed with the unique indefinite K\"ahler
  metric of index $2s$ such that the projection $M\to M/S^1$ becomes a semi-Riemannian submersion 
  (see \cite{bar}).
  $\mathbb CH^m_s$ is called the {\it complex pseudo-hyperbolic space}.
  Notice that $\mathbb CH^m_s$ has constant holomorphic sectional curvature $-4$
  and the curvature tensor is given by $R(X,Y,X,Y)=-g(X,X)g(Y,Y)+g(X,Y)^2-3g(I_0X,Y)^2$, where
  $I_0$ is the natural complex structure on $\mathbb CH^m_s$. 
   $\mathbb CH^m_s$ is a homogeneous space, namely we have   (see \cite{wol})
   $\mathbb CH^m_s=SU(s+1,m-s)/S(U(1)U(s,m-s))$ and\\
   $\mathbb CH^{2m+1}_{2s+1}=Sp(s+1,m-s)/U(1)Sp(s,m-s)$.
\end{defn}

We denote by $H^n(-4)$ the 
hyperbolic space with sectional curvature $-4$,
by $\mathbb HH^n$ the quaternionic hyperbolic 
space of real dimension $4n$ with 
quaternionic sectional curvature $-4$.

Many explicit examples of semi-Riemannian submersions with totally geodesic fibres can be
given following a standard construction (see \cite{bes} for Riemannian case).
Let $G$ be a Lie group and $K$, $H$ two compact Lie subgroups of
$G$ with $K\subset H$. Let $\pi:G/K\to G/H$ be the associated bundle 
with fibre $H/K$ to the $H-$principal 
bundle $p:G\to G/H$. Let $\tt g$ be the Lie algebra
of $G$ and $\tt k\subset\tt h$ the corresponding Lie subalgebras of $K$ and $H$.
We choose an $Ad(H)-$invariant complement $\tt m$ to $\tt h$ in $\tt g$,
and an $Ad(K)-$invariant complement $\tt p$ to $\tt k$ in $\tt h$.
An $ad(H)$-invariant nondegenerate bilinear symmetric form on
$\tt m$ defines a $G$-invariant semi-Riemannian metric $g'$ on $G/H$
and an $ad(K)$-invariant nondegenerate bilinear symmetric form on
$\tt p$ defines a $H$-invariant semi-Riemannian metric $\hat g$ on $H/K$.
The orthogonal direct sum for these nondegenerate bilinear symmetric
forms on $\tt p\oplus\tt m$ defines a $G$-invariant semi-Riemannian
metric $g$ on $G/K$.
The following theorem is proved in \cite{bes}.
\bth 
 The map $\pi:(G/K,g)\to(G/H,g')$ is a semi-Riemannian submersion 
with totally geodesic fibres.
\eth
Using this theorem we get the following examples.
\begin{exmp}
   Let $G=SU(1,n)$, $H=S(U(1)U(n))$, $K=SU(n)$.
   We have the semi-Riemannian submersion
   $$H^{2n+1}_1=SU(1,n)/SU(n)\to\mathbb CH^n=SU(1,n)/S(U(1)U(n)).$$
\end{exmp}
\begin{exmp}
   Let $G=Sp(1,n)$, $H=Sp(1)Sp(n)$, $K=Sp(n)$.
   We get the semi-Riemannian submersion
   $$H^{4n+3}_3=Sp(1,n)/Sp(n)\to\mathbb HH^n=Sp(1,n)/Sp(1)Sp(n).$$
\end{exmp}
\begin{exmp}
   Let $G=Spin(1,8)$, $H=Spin(8)$, $K=Spin(7)$.
   We have the semi-Riemannian submersion
   $$H^{15}_7=Spin(1,8)/Spin(7)\to H^8(-4)=Spin(1,8)/Spin(8).$$
\end{exmp}
\begin{exmp}
    Let $G=Sp(1,n)$, $H=Sp(1)Sp(n)$, $K=U(1)Sp(n)$.
   We obtain the semi-Riemannian submersion
   $$\mathbb CH^{2n+1}_1=Sp(1,n)/U(1)Sp(n)\to\mathbb HH^n=Sp(1,n)/Sp(1)Sp(n).$$
\end{exmp}

\begin{defn}
     Two semi-Riemannian submersions $\pi,\pi':(M,g)\to (B,g')$
     are called {\it equivalent} if there is an
     isometry $f$ of $M$ which induces an isometry $\tilde f$ of $B$ so that
     $\pi'\circ f=\tilde f\circ\pi$. In this case the pair $(f,\tilde f)$ is called a
     {\it bundle isometry}.
\end{defn}

We shall need the following theorem, which is the semi-Riemannian version of 
theorem 2.2 in \cite{esc}.
\bth\label{t:esc}
     Let $\pi_1,\pi_2:M\to B$ be semi-Riemannian submersions from a connected
     complete semi-Riemannian manifold onto a semi-Riemannian manifold.
     Assume the fibres of these submersions are connected and totally 
     geodesic. Suppose $f$ is an isometry of $M$ which satisfies the following
     two properties at a given point $p\in M$:
     \begin{itemize}
          \item[(1)] $f_{*p}:T_pM\to T_{f(p)}M$ maps $\mathcal{H}_{1p}$ onto
          $\mathcal{H}_{2f(p)}$, where  $\mathcal{H}_{i}$ denotes the horizontal
          distribution of $\pi_i$, $i\in\{1,2\}$;
          \item[(2)] For every $E$, $F\in T_pM$, $f_*A_{1E}F=A_{2{f_*E}}f_*F$,
          where $A_i$ are the integrability tensors associated with $\pi_i$.
     \end{itemize}
     Then $f$ induces an isometry $\Tilde{f}$ of $B$ so that the pair
     $(f,\Tilde{f})$ is a bundle isometry between $\pi_1$ and $\pi_2$.
     In particular, $\pi_1$ and $\pi_2$ are equivalent.
\eth

\section{Semi-Riemannian submersions with totally geodesic fibres}

\bpr\label{curbneg}
    If $\pi:H^m_s\to B^n$ is a semi-Riemannian submersion with totally
    geodesic fibres from an $m$-dimensional pseudo-hyperbolic
    space of index $s$ onto an $n$-dimensional 
    Riemannian manifold then $m=n+s$,
    the induced metrics on fibres are negative definite
    and $B$ has negative sectional curvature.
\epr
\begin{proof}
     By lemma \ref{lemtrei}-b), we get
     $g(A_XV,A_XV)=-g(X,X)g(V,V)\geq 0$ for every
     horizontal vector $X$ and for every vertical vector $V$. 
     Therefore $g(V,V)\leq 0$ for every vertical vector $V$.
     By lemma \ref{lemtrei}-c), we have 
     $R'(\pi_*X,\pi_*Y,\pi_*X,\pi_*Y)=-g'(\pi_*X,\pi_*X)g'(\pi_*Y,\pi_*Y)
     +g'(\pi_*X,\pi_*Y)^2+3g(A_XY,A_XY)<0$ 
     for every  linearly independent horizontal
     vectors $X$ and $Y$.
\end{proof}
\bpr\label{simply}
   Let $\pi:(M^{n+s}_s,g)\to (B^n,g')$ be a semi-Riemannian submersion
   from an $(n+s)$-dimensional semi-Riemannian manifold of index $s\geq 1$ onto
   an $n$-dimensional Riemannian manifold.
   We suppose $M$ is geodesically complete and simply connected.
   Then $B$ is complete and simply connected.
   If moreover $B$ has nonpositive curvature then the fibres are 
   simply connected.
\epr
\begin{proof}
    Since  $M$ is geodesically complete, the base space $B$ is complete.\\
    Let $\tilde g$ be the Riemannian metric on $M$ defined by
       $$\tilde g(E,F)=g(hE,hF)-g(vE,vF)$$
     for every $E$, $F$ vector fields on $M$.
    Since $\tilde g$ is a horizontally complete Riemannian metric
    (this means that any maximal horizontal geodesic 
    is defined on the entire real line) and 
    $B$ is a complete Riemannian manifold then $\mathcal H$ is an 
    Ehresmann connection for $\pi$ 
    (see theorem 1 in \cite{zhu}).
    By theorem 9.40 in \cite{bes}, it 
    follows $\pi:M\to B$ is a locally
    trivial fibration and we have an exact homotopy sequence
    $$\cdot\cdot\cdot\to\pi_2(B)\to\pi_1(fibre)\to\pi_1(M)\to\pi_1(B)\to 0$$
    Since $M$ is simply connected, we have $\pi_1(B)=0$.\\
    If $B$ has nonpositive curvature, then  $\pi_2(B)=0$
    by theorem of Hadamard.
    It follows $\pi_1(fibre)=0$.
\end{proof}
\bth\label{ajth}
  If $\pi:H^m_s\to B^n$ is a semi-Riemannian submersion with totally
  geodesic fibres from a pseudo-hyperbolic
  space of index $s>1$ onto a Riemannian manifold 
  then $B$ is a Riemannian symmetric space of
  rank one, noncompact and simply connected, 
  any fibre is diffeomorphic to $S^s$ and  $s\in\{3,7\}$.
\eth
\begin{proof}
     In order to prove that  $B$ is a locally symmetric space we need to check that
     $\na 'R'\equiv 0$.

     Let $X_0'$, $X'$, $Y'$, $Z'$ be vector fields on $B$ and let 
     $X_0$, $X$, $Y$, $Z$ be the horizontal liftings of $X_0'$, $X'$, $Y'$, $Z'$
     respectively. 
     By definition of the covariant derivative we have
\bea\label{rel*}
    (\na'_{X'_0}R')(X',Y',Z')&=&\na'_{X'_0}R'(X',Y')Z'-R'(\na'_{X'_0}X',Y')Z'\\
           &&-R'(X',\na'_{X'_0}Y')Z'-R'(X',Y')\na'_{X'_0}Z'.\no
\eea
In order to prove that the curvature tensor $R'$ of the base space is parallel,
we have to lift all vector fields in relation \eqref{rel*}.
By lemma \ref{lemdoi}, the horizontal liftings of
 $\na'_{X'_0}X'$, $\na'_{X'_0}Y'$
and $\na'_{X'_0}Z'$ are  $h\na_{X_0}X$, $h\na_{X_0}Y$ and $h\na_{X_0}Z$, 
respectively.
 
     We denote by $R^h(X,Y)Z$ the horizontal lifting of $R'(X',Y')Z'$. 
     The convention for Riemann tensor used here is 
     $R(X,Y)=\na_X\na_Y-\na_Y\na_X-\na_{[X,Y]}$.
     O'Neill's equation $vi)$ gives us the
     following relation
          $$R^h(X,Y)Z=h(R(X,Y)Z)+2A_ZA_XY-A_XA_YZ-A_YA_ZX.$$
     Using this relation we compute
\bea\label{rel***}
      (\na'_{X_0}R')(X',Y',Z')
         &=& \pi_*[ h\na_{X_0}(R^h(X,Y)Z)-R^h(h\na_{X_0}X,Y)Z\no\\
         && -R^h(X,h\na_{X_0}Y)Z-R^h(X,Y)h\na_{X_0}Z]\no\\
         &=& \pi_*[h\na_{X_0}h(R(X,Y)Z)-hR(h\na_{X_0}X,Y)Z\\
         && -hR(X,h\na_{X_0}Y)Z-hR(X,Y)h\na_{X_0}Z\no\\
         && +2(h\na_{X_0}A_ZA_XY-A_{h\na_{X_0}Z}A_XY\no\\
         && -A_ZA_{h\na_{X_0}X}Y-A_ZA_Xh\na_{X_0}Y)\no\\
         && -(h\na_{X_0}A_XA_YZ-A_{h\na_{X_0}X}A_YZ\no\\
         && -A_XA_{h\na_{X_0}Y}Z-A_XA_Yh\na_{X_0}Z)\no\\
         && -(h\na_{X_0}A_YA_ZX-A_{h\na_{X_0}Y}A_ZX\no\\
         && -A_YA_{h\na_{X_0}Z}X-A_YA_Zh\na_{X_0}X)].\no
\eea

   Since $H^m_s$ has constant curvature, we have $R(X,Y,Z,U)=0$ for every
   vertical vector $U$ and for every horizontal vector fields $X$, $Y$, 
   $Z$. This implies
   $R(X,Y)Z$ is horizontal and 
   $R(X,U)Y$, $R(U,X)Y$, $R(X,Y)U$ are vertical. Hence\\
   $\pi_*(\na_{X_0}h(R(X,Y)Z)-hR(h\na_{X_0}X,Y)Z-hR(X,h\na_{X_0}Y)Z
            -hR(X,Y)h\na_{X_0}Z)
   =\pi_*(\na_{X_0}R(X,Y)Z)- \pi_*(R(\na_{X_0}X,Y)Z-R(v\na_{X_0}X,Y)Z)
   -\pi_*(R(X,Y)\na_{X_0}Z-R(X,Y)v\na_{X_0}Z)=\pi_*[(\na_{X_0}R)(X,Y,Z)]$.
   Since $H^m_s$ has constant curvature, we get 
   $(\na_{X_0}R)(X,Y,Z)=0$. 
   So the sum of the first four terms in relation \eqref{rel***}
   is zero.

     We have
     $h\na_{X_0}A_ZA_XY-A_{h\na_{X_0}Z}A_XY-A_ZA_{h\na_{X_0}X}Y-A_ZA_Xh\na_{X_0}Y=$\\
     $h((\na_{X_0}A)_Z(A_XY))-A_Z(v(\na_{X_0}A)_XY)$.\\
     For the case of totally geodesic fibres, O'Neill's equation $v)$ becomes
     $$R(X,Y,Z,U)=g((\na_ZA)_XY,U).$$
     By lemma \ref{l:1} (f) and by
     the hypothesis of constant curvature total space we get
     $$g((\na_ZA)_XU,Y)=g((\na_ZA)_XY,U)=0$$
     for every horizontal vector fields $X$, $Y$, $Z$
     and for every vertical vector field $U$.
     It follows $h(\na_ZA)_XU=0$ and $v(\na_ZA)_XY=0$
     for every horizontal vector fields $X$, $Y$, $Z$
     and for every vertical vector field $U$.
     Therefore $h((\na_{X_0}A)_Z(A_XY))=0$ and $v((\na_{X_0}A)_XY)=0$.
     This implies 
       $$h\na_{X_0}A_ZA_XY-A_{h\na_{X_0}Z}A_XY-A_ZA_{h\na_{X_0}X}Y-A_ZA_Xh\na_{X_0}Y=0.$$
     By circular permutations of $(X,Y,Z)$ in the last relation we get
     $$h\na_{X_0}A_XA_YZ-A_{h\na_{X_0}X}A_YZ-A_XA_{h\na_{X_0}Y}Z-A_XA_Yh\na_{X_0}Z=0,$$
     $$h\na_{X_0}A_YA_ZX-A_{h\na_{X_0}Y}A_ZX-A_YA_{h\na_{X_0}Z}X-A_YA_Zh\na_{X_0}X=0.$$
     So the sum of all terms in relation \eqref{rel***} is zero.\\
     We proved that $(\na '_{X'_0}R')(X',Y',Z')=0$ for every vector fields $X_0'$, $X'$, $Y'$, $Z'$,
     so $B$ is a locally symmetric space.
     By proposition \ref{simply}, $B$ is simply connected and complete.
     Therefore $B$ is a Riemannian symmetric space.
     By proposition \ref{curbneg}, $B$ has negative sectional curvature. Hence $B$ 
     is a noncompact Riemannian symmetric space of rank one.

   Let $b\in B$. Since $\pi^{-1}(b)$ is a totally geodesic submanifold
   of a geodesically complete manifold, $\pi^{-1}(b)$
   is itself geodesically complete.
   Since
   $R'(X',Y',X',Y')\leq 0$ for every $X'$, $Y'$ tangent vectors to $B$,
   we have $\pi_1(fibre)=0$, by proposition \ref{simply}.
   Since $(\pi^{-1}(b),\hat g)$ is a complete, simply connected semi-Riemannian
   manifold of dimension $r$ and of index $r$ and with constant sectional curvature $-1$,
   it follows $(\pi^{-1}(b),\hat g)$ is isometric to $H^s_s$ 
   (see Proposition 23 from page 227 in \cite{onei}).
    Hence  any fibre is diffeomorphic to $S^s$.

   We shall prove below that the tangent bundle of any fibre is trivial.
   {From} a well known result of Adams it follows that $s\in\{1,3,7\}$.
\ble\label{ajunu}
   The tangent bundle of any fibre is trivial.
\ele
\begin{proof}[Proof of lemma \ref{ajunu}]
   Since $g(A_XV,A_XV)=-g(X,X)g(V,V)$, we have that\\
   $A_X:\mathcal V\to\mathcal H$, $V\mapsto A_XV$ is an injective map and
   $\dim \mathcal V\leq\dim \mathcal H$, if $g(X,X)\not=0$.\\
   For any  horizontal vector field $X$, 
   we denote by
   $A^*_X:\mathcal H\to\mathcal V$ the map given by $A^*_X(Y)=A_XY$.
   By O'Neill's equation $iv)$, we have $g(A_XV,A_XW)=-g(X,X)g(V,W)$ for every 
   vertical vector fields $V$ and $W$. Hence, by lemma \ref{l:1} $(d)$, we get
   $A^*_XA_XV=g(X,X)V$ for every vertical vector field $V$. 
   If $g(X,X)\not=0$ anywhere then $A^*_X$ is surjective
   and hence
   $\dim\mathcal V=
    \dim\mathcal H-\dim\ker A^*_X$. 
   By lemma \ref{l:1} $(d)$, we have $A_XX=0$. This implies
	$\dim\ker A^*_X\geq 1$.\\
   Let $b\in B$ and $x\in T_bB$ with $g(x,x)=1$.
   We denote by $X$ the horizontal lifting along the fibre $\pi^{-1}(b)$
   of the vector $x$.
   Let $p$ an arbitrary point in $\pi^{-1}(b)$ and
   let $\{X(p),y_1,\cdot\cdot\cdot,y_l\}$ be
   an orthonormal basis of the vector space $\ker A^*_{X(p)}$.
   Since $\pi_{*p}$ sends isometrically $\mathcal H_p$ into $T_bB$ we have 
   $\{\pi_*X(p),\pi_*y_1,\cdot\cdot\cdot ,\pi_*y_l\}$ is a linearly independent
   system which can be completed to a basis of $T_bB$ with a system of
   vectors $\{x_{l+1},\cdot\cdot\cdot,x_{n-1}\}$.
   Let $X,X_1,X_2,\cdot\cdot\cdot,X_{n-1}$ be the horizontal liftings
   along the fibre $\pi^{-1}(b)$ of
   $x=\pi_*X(p),\pi_*y_1,\cdot\cdot\cdot ,\pi_*y_l,
                 x_{l+1},\cdot\cdot\cdot,x_{n-1}$ respectively.

   By lemma \ref{lemtrei}, we have for every $q\in\pi^{-1}(b)$ and 
   for every $i\in\{1,\dots,l\}$
\bea
   3g(A_{X(q)}X_i(q),A_{X(q)}X_i(q))&=& R'(\pi_*X(q),\pi_*X_i(q),\pi_*X(q),\pi_*X_i(q))\no\\
   &&-R(X(q),X_i(q),X(q),X_i(q))\no\\
   &=& R'(x,\pi_*y_i,x,\pi_*y_i)\no\\
    &&+g(X(q),X(q))g(X_i(q),X_i(q))-g(X(q),X_i(q))^2\no\\
   &=&R'(x,\pi_*y_i,x,\pi_*y_i)\no\\
   &&+g'(\pi_*X(q),\pi_*X(q))g'(\pi_*X_i(q),\pi_*X_i(q))\no\\
   &&-g'(\pi_*X(q),\pi_*X_i(q))^2\no\\
   &=&3g(A_{X(p)}X_i(p),A_{X(p)}X_i(p))\no\\
   &=&0.\no
\eea
   Since the induced metrics on fibre $\pi^{-1}(b)$ are negative definite, we get
   $A_{X(q)}X_i(q)=0$.\\
   By lemma \ref{l:1} (a), we have $A_{X(q)}{X(q)}=0$.
   We proved that $\{X(q),X_1(q),\dots,X_l(q)\}\subset \ker A^*_{X(q)}$.
   Since $\pi_{*q}$ sends isometrically $\mathcal H_q$ into $T_bB$, we get
   $\{X(q),X_1(q),\dots,X_l(q)\}$
   is a basis of the vector space 
   $\ker A^*_{X(q)}$ for every point
   $q\in\pi^{-1}(b)$.\\
   Let $V_{l+1}=A^*_XX_{l+1},\cdot\cdot\cdot,V_{n-1}=A^*_XX_{n-1}$ 
   be tangent vector fields to the fibre $\pi^{-1}(b)$.
   We denote by $Q_q$ the vector subspace of $\mathcal H_q$ spanned by
   $\{X_{l+1}(q),X_{l+2}(q),\cdot\cdot\cdot,X_{n-1}(q)\}$.
   Let $\tilde g$ be the Riemannian metric on $\pi^{-1}(b)$ given by
   $\tilde g(V,W)=-g(V,W)$ for every $V$, $W$ vector fields  tangent to 
   $\pi^{-1}(b)$.
   Since $\dim \mathcal V_q = \dim Q_q$ and 
   $g(A_XV,A_XV)=\tilde g(V,V)$,
   we get
     $A_{X(q)}:(\mathcal V_q,\tilde g)\to(Q_q,g)$
   is an isometry for every $q\in\pi^{-1}(b)$.\\ 
   So $\{V_{l+1},\cdot\cdot\cdot,V_{n-1}\}$ 
   is a global frame for the tangent bundle of
   $\pi^{-1}(b)$. 
   It follows the tangent bundle of the fibre $\pi^{-1}(b)$
   is trivial.
\end{proof}

This ends the proof of theorem \ref{ajth}.
\end{proof}

By the classification of the Riemannian
symmetric spaces of rank one of noncompact type, we
have $B$ is isometric to one of the following spaces:
\begin{itemize}
\item[1)] $H^n(c)$ real hyperbolic space with constant sectional curvature $c$;
\item[2)] $\mathbb CH^k(c)$ complex hyperbolic space
           with holomorphic sectional curvature $c$;
\item[3)] $\mathbb HH^k(c)$ quaternionic hyperbolic space
          with quaternionic sectional curvature $c$;
\item[4)] $\mathbb CaH^2(c)$ Cayley hyperbolic plane
          with Cayley sectional curvature $c$.
\end{itemize}

This will give us more information about the relation between
the dimension of fibres and the geometry of base space.
\bpr\label{geo}
   Let $\pi:H^{n+s}_s\to B^n$ be a semi-Riemannian submersion with totally
   geodesic fibres.
\begin{itemize}
\item[$a)$] If $s=3$ then  $n=4k$ and $B^n$ is isometric to $\mathbb HH^k$.
\item[$b)$] If $s=7$ then we have one of the following situations:
           \begin{itemize}
   	      \item[i)] $n=8$ and $B^n$ is isometric to 
                 $H^8(-4)$; or
              \item[ii)] $n=16$ and $B^n$ is isometric to
                $\mathbb CaH^2$.
           \end{itemize}
\end{itemize}
\epr
\begin{proof}
Let $Y$, $Z$ be two linear independent horizontal vectors and 
let $Y'=\pi_*Y$,  $Z'=\pi_*Z$.
By proposition \ref{curbneg}, 
the metric induced on fibres are negative definite.
This implies $g(A_ZY,A_ZY)\leq 0$. By lemma \ref{lemtrei}, we get
$$K'(Z',Y')=\f{R'(Z',Y',Z',Y')}{g'(Z',Z')g'(Y',Y')-g'(Z',Y')^2}=
-1+\f{3g(A_ZY,A_ZY)}{g(Z,Z)g(Y,Y)-g(Z,Y)^2}\leq -1.$$
By Schwartz inequality applied to the positive definite scalar product 
induced on $\mathcal H$, we have
$$-g(A_ZY,A_ZY)=g(A_ZA_ZY,Y)\leq\sqrt{g(A_ZA_ZY,A_ZA_ZY)}\sqrt{g(Y,Y)}.$$
By lemma \ref{lemtrei}, we get 
$-g(A_ZY,A_ZY)\leq\sqrt{-g(A_ZY,A_ZY)g(Z,Z)}\sqrt{g(Y,Y)}.$ Thus
$-g(A_ZY,A_ZY)\leq g(Z,Z)g(Y,Y)$.
Therefore $K'(Z',Y')=-1+\f{3g(A_ZY,A_ZY)}{g(Z,Z)g(Y,Y)}\geq -4$
for every  orthogonal vectors $Z'$ and $Y'$. 

We proved that $-4\leq K'\leq -1$.\\
We shall prove that 
if the base space $B$ has constant curvature $c$ then $c=-4$.
It is sufficient to see that for any point $b\in B$
there is a $2-$plane $\al\in T_bB$ such that $K(\al)=-4$.
We choose $\al=\{\pi_*Z,\pi_*A_ZV\}$ where $Z$ is a horizontal vector and 
$V$ is a vertical vector.
By lemma \ref{lemtrei}, we have 
\bea\label{eq**}
  R'(\pi_*Z,\pi_*A_ZV,\pi_*Z,\pi_*A_ZV)&=&R(Z,A_ZV,Z,A_ZV)\\
       && +3g(A_Z(A_ZV),A_Z(A_ZV)).\no
\eea
We notice that $Z$ and $A_ZV$ are orthogonal, because, by lemma \ref{l:1}, we have 
$g(Z,A_ZV)=-g(A_ZZ,V)=0$.\\
By lemma \ref{lemtrei}, we have 
$$g(A_XU,A_XU)=-g(X,X)g(U,U)$$
for every horizontal vector $X$ and for every vertical vector $U$.
By lemma \ref{l:1} (d), we get $g(A_XA_XU,U)=g(X,X)g(U,U)$. Hence, by polarization,
we find $A_XA_XU=g(X,X)U$ for every horizontal vector $X$ and for every vertical vector $U$.
Therefore the relation \eqref{eq**} becomes 
\bea
  R'(\pi_*Z,\pi_*A_ZV,\pi_*Z,\pi_*A_ZV)&=&-g(Z,Z)g(A_ZV,A_ZV)+3g(Z,Z)^2g(V,V)\no\\
&=&4g(Z,Z)^2g(V,V)\no\\
&=&-4(g'(\pi_*Z,\pi_*Z)g'(\pi_*A_ZV,\pi_*A_ZV)-g'(\pi_*Z,\pi_*A_ZV)^2).\no
\eea
Then
$K'(\pi_*Z,\pi_*A_ZV)=-4$. 
Therefore if the base space $B$ has constant curvature $c$ then
$c=-4$.

Let $X$ be a horizontal vector field. By lemma \ref{lemtrei},
     $Y\in\ker A^*_X$ if and only if
    $$R'(\pi_*X,\pi_*Y,\pi_*X,\pi_*Y)=-g'(\pi_*X,\pi_*X)g'(\pi_*Y,\pi_*Y)+
    g'(\pi_*X,\pi_*Y)^2.$$
    For every $X'\in T_{\pi(p)}B$, we denote by
      $$\mathcal L_{X'}=\{Y'\in T_{\pi(p)}B\ |
    \  R'(X',Y',X',Y')=-g'(X',X')g'(Y',Y')+
        g'(X',Y')^2\} .$$
    With this notation, $\pi_*(\ker A^*_{X(p)})=\mathcal L_{\pi_*X(p)}$.
    Since $\pi_*$ sends isometrically $\mathcal H_p$ into $T_{\pi(p)}B$, 
    we have $\dim \mathcal H-\dim \mathcal V= \dim\ker 
    A^*_{X(p)}=\dim \mathcal L_{\pi_*X(p)}$\\
    We compute $\dim \mathcal L_{X'}$ from the geometry of $B$. We have
     the following possibilities for $B$:
\begin{case} $B=H^k(-4)$.\\
        The curvature tensor of hyperbolic space $H^k(-4)$ is given by
         $$R'(X',Y',X',Y')=-4(g'(X',X')g'(Y',Y')-g'(X',Y')^2).$$
        We have $\mathcal L_{X'}=\{\la X'\ |\ \la\in\mathbb R\}$. Hence
        $\dim \mathcal L_{X'}=1$. It follows 
         $\dim \mathcal H=\dim \mathcal V+1$.

       If $s=3$ then $B^n$ is isometric to $H^4(-4)$, which falls in
       the case a), since $H^4(-4)$ is isometric to $\mathbb HH^1$.

       If $s=7$ then $\dim \mathcal H=8$ and this is the case b ii).
\end{case}
\begin{case} $B=\mathbb CH^k$.\\
     Let $I_0$ be the natural complex structure on $\mathbb CH^k$.
     The curvature tensor of complex hyperbolic space
     $\mathbb CH^k$ with $-4\leq K'\leq -1$ is given by 
      %(see \cite{kob})
         $$R'(X',Y',X',Y')=-(g'(X',X')g'(Y',Y')-g'(X',Y')^2+3g'(I_0X',Y')^2).$$
     We get $\mathcal L_{X'}=\{I_0X'\}^{\perp}$. So
     $\dim\mathcal L_{X'}=2k-1=\dim\mathcal H-1$. It follows 
         $\dim\mathcal V=1$.
\end{case}
\begin{case} $B=\mathbb HH^k$.\\
      Let $\{I_0,J_0,K_0\}$ be local almost complex structures
     which rise to the quaternionic structure on $\mathbb HH^k$.
     The curvature tensor of the quaternionic hyperbolic space
     $\mathbb HH^k$ with $-4\leq K'\leq -1$ (see \cite{ish}) is given by 
\bea     
    R'(X',Y',X',Y')&=&
   -g'(X',X')g'(Y',Y')+g'(X',Y')^2\no\\
   &&-3g'(I_0X',Y')^2-3g'(J_0X',Y')^2-3g'(K_0X',Y')^2.\no
\eea
     It follows that $Y'\in\mathcal L_{X'}$ if and only if 
     $g'(I_0X',Y')=g'(J_0X',Y')=g'(K_0X',Y')=0$. Therefore
      $\mathcal L_{X'}=\{I_0X',J_0X',K_0X'\}^{\perp}$. Hence
     $\dim\mathcal L_{X'}=4k-3=\dim\mathcal H-3$. We get
         $\dim\mathcal V=3$.
\end{case}
\begin{case}
     $B=\mathbb CaH^2$.\\
     Let $\{I_0,J_0,K_0,M_0,M_0I_0,M_0J_0,M_0K_0\}$ be 
     local almost complex structures
     which rise to the Cayley structure
     on $\mathbb CaH^2$ Cayley hyperbolic plane. 
     The curvature tensor of the Cayley plane $\mathbb CaH^2$ with $-4\leq K'\leq -1$
     (see \cite{cou}) is given by 
\bea
   R'(X',Y',X',Y')&=&
   -g'(X',X')g'(Y',Y')+g'(X',Y')^2-3g'(I_0X',Y')^2\no\\
   &&-3g'(J_0X',Y')^2-3g'(K_0X',Y')^2-3g'(M_0I_0X',Y')^2\no\\
   &&-3g'(M_0J_0X',Y')^2-3g'(M_0K_0X',Y')^2.\no
\eea
   We get
$$\mathcal L_{X'}=\{I_0X',J_0X',K_0X',
M_0X',M_0I_0X',M_0J_0X',M_0K_0X'\}^{\perp}.$$
    So  $\dim \mathcal L_{X'}=\dim\mathcal H-7$. It follows
    $\dim \mathcal V=7$.
\end{case}
\end{proof}

Summarizing all of the above, we obtain our main classification result.

\bmth\label{main}
    Let $\pi:H^m_s\to B$ be a semi-Riemannian submersion with totally
    geodesic fibres from a pseudo-hyperbolic space
    onto a Riemannian manifold. Then the semi-Riemannian submersion $\pi$ is equivalent
    to one of the following canonical semi-Riemannian submersions, 
    given by examples 1)-3):
\begin{itemize}
	\item[$(a)$] $H^{2k+1}_1\to\mathbb CH^k$,
	\item[$(b)$] $H^{4k+3}_3\to\mathbb HH^k$,
	\item[$(c)$] $H^{15}_7\to H^8(-4)$.
\end{itemize}
\emth
\begin{proof}
    The index of the pseudo-hyperbolic space cannot be $s=0$. 
    Indeed, by lemma \ref{lemtrei}, for $s=0$,
    we get $0\leq g(A_XV,A_XV)=-g(X,X)g(V,V)\leq 0$ 
    for every  horizontal vector $X$ and 
    for every  vertical  vector $V$. But this is not possible.

    By \cite{mag}, any semi-Riemannian submersion with
    totally geodesic fibres, from 
    a pseudo-hyperbolic space of index $1$ onto a Riemannian manifold 
    is equivalent to the canonical 
    semi-Riemannian submersion $H^{2k+1}_1\to\mathbb CH^k$.

    It remains to study the case $s>1$.
    By theorem \ref{ajth} and proposition \ref{geo},
    any semi-Riemannian
    submersion with totally geodesic fibres
    from a pseudo-hyperbolic space of index $s>1$
    onto a Riemannian manifold is one of the following types: 
\begin{itemize}
\item[(1)]  $H^{4k+3}_3\to\mathbb HH^k$, or
\item[(2)]  $H^{15}_7\to H^8(-4)$, or
\item[(3)]  $H^{23}_7\to\mathbb CaH^2$.
\end{itemize}
    
     In order to prove that any two semi-Riemannian submersions in one of the
     categories (1) or (2) are equivalent we 
     shall modify Ranjan's argument (see \cite{ran})
     to our situation. In the category (3), 
     we shall prove there are no such semi-Riemannian submersions
     with totally geodesic fibres.

     First, we shall prove the uniqueness in the case $H^{4k+3}_3\to\mathbb HH^k$.
     Let $p\in H^{4k+3}_3$ and let $\mathcal U:\mathcal V_p\to End(\mathcal H_p)$
     the map given by $\mathcal U(v)(x)=A_xv$ for every $v\in\mathcal V_p$ and
     for every $x\in\mathcal H_p$.
     We denote $\mathcal U(v)$ by $A^v$. It is trivial to see that $A^v$ is 
     skew-symmetric (i.e. $g(A^vx,y)=-g(x,A^vy)$).
     The O'Neill's equation $g(A_xv,A_xv)=-g(x,x)g(v,v)$ becomes
     $g(A^vx,A^vx)=-g(x,x)g(v,v)$. This implies $g(A^vA^vx,x)=g(x,x)g(v,v)$. 
     Hence, by polarization in $x$, we have $g(A^vA^vx,y)=g(x,y)g(v,v)$ for every
     $y\in\mathcal H_p$.
     So $A^vA^vx=g(v,v)x$. Again by polarization we get $A^vA^w+A^wA^v=2g(v,w)Id$.
     Let $\tilde g$ be the Riemannian metric given by $\tilde g(v,w)=-g(v,w)$ for
     every $v$, $w\in\mathcal V_p$.
     It follows $A^vA^w+A^wA^v=-2\tilde g(v,w)Id_{\mathcal V_p}$.
     This is the condition which
     allows us to extend $\mathcal U$ to a representation of the Clifford algebra
     $Cl(\mathcal V_p,\tilde g_p)$ of $\mathcal V_p$. We also denote by $\mathcal U$
     the extension of $\mathcal U$.
     Since $\dim\mathcal V_p=3$ and $\tilde g_p$ is positive definite,
     $Cl(\mathcal V_p,\tilde g_p)$ has at most two types
     of irreducible representations. We notice that $\mathcal H_p$ is a
     $Cl(\mathcal V_p,\tilde g_p)$-module which splits in simple modules of dimension
     $4$. The next step is to show that any two such simple modules in decomposition
     of $\mathcal H_p$ are equivalent. Let $\{v_1,v_2,v_3\}$ be an orthonormal basis 
     of $(\mathcal V_p,\tilde g_p)$.
     Since the affiliation of a simple $Cl(\mathcal V_p,\tilde g_p)$-module
     to one of the two possible types is decided by the action of $v_1v_2v_3$, it is
     sufficient to check that $A^{v_1}A^{v_2}A^{v_3}=Id_{\mathcal V_p}$.

      Consider the function $x\mapsto g(A^{v_1}A^{v_2}A^{v_3}x,x)$ defined on the unit
      sphere in $\mathcal H_p$. We have
      $g(A^{v_1}A^{v_2}A^{v_3}x,x)=
            -g(A^{v_2}A^{v_3}x,A^{v_1}x)=g(A_xA_{A_x{v_3}}v_2,v_1)$. 
      A straightforward computation shows that $A_xA_{A_x{v_3}}v_2$ is orthogonal to
      $v_2$ and $v_3$. Hence $A_xA_{A_x{v_3}}v_2$ is a multiple of $v_1$.

      By polarization of the relation
      $A_xA_xv=g(x,x)v$, we get $A_xA_y+A_yA_x=2g(x,y)Id$
      for every horizontal vectors $x$ and $y$ .
      In particular, we have 
      $A_xA_{A_x{v_3}}v_2=-A_{A_x{v_3}}A_xv_2+2g(x,A_xv_3)v_2=-A_{A_x{v_3}}A_xv_2$.

      Let $S$ be the vector subspace of $\mathcal H_p$ spanned by 
      $\{x, A_xv_1, A_xv_2, A_xv_3\}$. By lemma \ref{lemtrei}, we get 
      $K'(\pi_*x,\pi_*A_xv_i)=-4$ for all $i\in\{1,2,3\}$.
      By geometry of $\mathbb HH^n$, there exists a 
      unique totally geodesic hyperbolic line $\mathbb HH^1$ passing through
      $\pi(p)$ such that $T_{\pi(p)}\mathbb HH^1=\pi_*S$. Notice that for every orthonormal 
      vectors $y,z\in T_{\pi(p)}\mathbb HH^1$, $K'(y,z)=-4$. In particular
      we have $K'(\pi_*A_xv_2,\pi_*A_xv_3)=-4$.
      Hence
        $g(A_{A_x{v_3}}A_xv_2,A_{A_x{v_3}}A_xv_2)=-1$. It follows that
        $A_xA_{A_x{v_3}}v_2=\pm v_1$. Hence $ g(A^{v_1}A^{v_2}A^{v_3}x,x)=\pm 1$
         for all unit vectors $x$.
      Since the function $x\mapsto g(A^{v_1}A^{v_2}A^{v_3}x,x)$ defined on the unit
      sphere in $\mathcal H_p$ is continuous, we get either
\begin{itemize}
 \item[(i)] $g(A^{v_1}A^{v_2}A^{v_3}x,x)=1$  for any unit horizontal vector $x$, or 
 \item[(ii)]  $g(A^{v_1}A^{v_2}A^{v_3}x,x)=-1$ for any unit horizontal vector $x$.
\end{itemize}

We may assume the case (i) holds.\\
If the case (ii) is happen, we replace the orthonormal basis $\{v_1,v_2,v_3\}$ of 
$(\mathcal V_p,\tilde g_p)$ with the orthonormal basis $\{v_1,v_2,-v_3\}$. So for this new basis
we are in the case (i).

      Since $A^{v_1}A^{v_2}A^{v_3}$ is an isometry, we have 
$g(A^{v_1}A^{v_2}A^{v_3}x,A^{v_1}A^{v_2}A^{v_3}x)g(x,x)=g(x,x)^2=1=g(A^{v_1}A^{v_2}A^{v_3}x,x)^2$
for all unit horizontal vectors $x$.

      So the Schwartz inequality for the scalar product $g | _{\mathcal H_p}$
      $$g(A^{v_1}A^{v_2}A^{v_3}x,x)^2\leq
          g(A^{v_1}A^{v_2}A^{v_3}x,A^{v_1}A^{v_2}A^{v_3}x)g(x,x)$$ 
      becomes equality. It follows that $A^{v_1}A^{v_2}A^{v_3}x=\la x$ for some $\la$.
    Because $A^{v_1}A^{v_2}A^{v_3}$ is an isometry and we assumed the case (i),
    it follows $\la=1$.
    We proved that $A^{v_1}A^{v_2}A^{v_3}x=x$ for all unit horizontal vectors $x$.
    Obviously, $A^{v_1}A^{v_2}A^{v_3}x=x$ for all $x\in\mathcal H_p$.

     Let $\pi':H^{4k+3}_3\to\mathbb HH^k$ be another semi-Riemannian submersion
     with totally geodesic fibres.
     For an arbitrary chosen point $q\in H^{4k+3}_3$, we consider horizontal and
     vertical subspaces $\mathcal H'_q$ and $\mathcal V'_q$.
     Let $\{v'_1,v'_2,v'_3\}$  be an orthonormal basis in $\mathcal V'_q$ such that 
     $v'_1v'_2v'_3$ acts on $\mathcal H'_q$ as $Id$.
     Let $L_1:\mathcal V'_q\to\mathcal V_p$ be the isometry given by $L_1(v'_i)=v_i$
     for all $i\in\{1,2,3\}$ and let 
     $Cl(L_1):Cl(\mathcal V'_q)\to Cl(\mathcal V_p)$ be the extension of $L_1$ 
     to the Clifford algebras.
     The composition $\mathcal U\circ Cl(L_1):Cl(\mathcal V'_q)\to End(\mathcal H_p)$ makes
     $\mathcal H_p$ to be a $Cl(\mathcal V'_q)$-module of dimension $4k$.
     Let $\mathcal H_p=\mathcal H_1\oplus\dots\oplus\mathcal H_k$ and
     $\mathcal H'_q=\mathcal H'_1\oplus\dots\oplus\mathcal H'_k$ be the decomposition
     of $\mathcal H_p$ and $\mathcal H'_q$ in simple $Cl(\mathcal V'_q)$-modules, respectively.
     For each $i$ there is  $f_i:\mathcal H'_i\to\mathcal H_i$ 
     an equivalence of $Cl(\mathcal V'_q)$-modules, which after a rescaling by a constant number
     is an isometry which preserves the O'Neill's integrability tensors.
     Taking the direct sum of all these isometries, we obtain an isometry
     $L_2:\mathcal H'_q\to\mathcal H_p$ which preserves the  O'Neill's integrability tensors.
     Therefore $L=L_1\oplus L_2:T_qH^{4k+3}_3\to T_pH^{4k+3}_3$ is an isometry which
     maps $\mathcal H'_q$ onto $\mathcal H_p$ and $A'$ onto $A$.
     Since $H^{4k+3}_3$ is a simply connected complete symmetric space,
     there is an isometry $f:H^{4k+3}_3\to H^{4k+3}_3$
     such that $f(q)=p$ and $f_{*q}=L$ (see corollary 2.3.14 in \cite{wol}).
     Therefore, by theorem \ref{t:esc}, we get $\pi$ and $\pi'$ are equivalent.

     Now, we shall prove that any two semi-Riemannian submersions
     $\pi,\pi':H^{15}_7\to H^8(-4)$ with totally geodesic fibres are equivalent.
     The proof is analogous to the case $(1)$, but it is easier.

     Let $p,q\in H^{15}_7$ and let 
     $\mathcal H_p$, $\mathcal V_p$ be the horizontal and vertical subspaces in
     $T_pH^{15}_7$ for $\pi$,
     let $\mathcal H'_q$, $\mathcal V'_q$ be the horizontal and vertical subspaces in
     $T_qH^{15}_7$ for $\pi'$.
     Let $\{v_1,\dots, v_7\}$  be an orthonormal basis of  $(\mathcal V_p,\tilde g_p)$
     and $\{v'_1,\dots, v'_7\}$  be an orthonormal basis of  $(\mathcal V'_q,\tilde g_q)$
     such that $A^{v_1}A^{v_2}\dots A^{v_7}=Id$ and
     $A^{v'_1}A^{v'_2}\dots A^{v'_7}=Id$.
     Since $\dim\mathcal V_p=7$, the irreducible $Cl(\mathcal V_p,\tilde g_p)$-modules are
     $8$-dimensional. Since $\dim\mathcal H_p=8$, we get $\mathcal H_p$ is simple.
     Because  $A^{v_1}A^{v_2}\dots A^{v_7}=Id$ and $A^{v'_1}A^{v'_2}\dots A^{v'_7}=Id$
     we get $\mathcal H'_q$ and $\mathcal H_p$ are $Cl(\mathcal V'_q)$-modules equivalent.
     Analogously to the case (1), we can construct an isometry 
     $L=L_1\oplus L_2:T_qH^{15}_7\to T_pH^{15}_7$, which
     map $\mathcal H'_q$ onto $\mathcal H_p$ and $A'$ onto $A$. This produces
     an isometry $f:H^{15}_7\to H^{15}_7$
     such that $f(q)=p$ and $f_{*q}=L$ (see corollary 2.3.14 in \cite{wol}).
     Again by theorem \ref{t:esc}, we get $\pi$ and $\pi'$ are equivalent.

   Now, we prove that there are no $\pi:H^{23}_7\to\mathbb CaH^2$ semi-Riemannian submersions
   with totally geodesic fibres. The proof is analogous to that of Ranjan 
   (see proposition 5.1 in \cite{ran}).

$\mathcal H_p$ becomes a $Cl(\mathcal V_p)$-module 
by considering the extension
of the map 
$\mathcal U:\mathcal V_p\to\mathrm{End}(\mathcal H_p)$, $\mathcal U(V)(X)=A_XV$
to the Clifford algebra $Cl(\mathcal V_p)$.
Here $Cl(\mathcal V_p,\Tilde g_p)$ denotes the Clifford algebra 
of $(\mathcal V_p,\Tilde g_p)$, $\Tilde g(U,V)=-g(U,V)$ 
for every $U$, $V\in\mathcal V_p$. 
Since $\Tilde g_p$ is positive definite, we have
$Cl(\mathcal V_p)\simeq\mathbb R(8)\oplus\mathbb R(8)$. 
Hence, $\mathcal H_p$ splits into
two $8$-dimensional irreducible $Cl(\mathcal V_p)$-modules.
Since the induced metrics on
fibres are negative definite 
we get $\pi^{-1}(\mathbb CaH^1)$ 
is totally geodesic in $H^{23}_{7}$ and isometric
to $H^{15}_7$, by theorem 2.5 in \cite{esco}. Here $\mathbb CaH^1$ denotes the 
Cayley hyperbolic line through $\pi_*X$;
We choose $S$ be the horizontal space of the restricted submersion
$\Tilde\pi:H^{15}_7\to\mathbb CaH^1=H^8(-4)$.
So for every 
$X\in\mathcal H_p$, $g(X,X)\not=0$ we find
an irreducible $Cl(\mathcal V_p)$-submodule $S$ of 
$\mathcal H_p$ passing through $X$. Since $\dim\mathcal V_p\geq 4$, we get a contradiction.

\end{proof}

    R. Escobales \cite{esco} classified Riemannian 
    submersions from complex projective spaces under
    the assumption that the fibres are connected, complex, 
    totally geodesic submanifolds.
    Using the main theorem \ref{main}, we obtain a 
    classification of semi-Riemannian submersions 
    from a complex pseudo-hyperbolic space onto
    a Riemannian manifold under
    the assumption that the fibres are connected, complex, 
    totally geodesic submanifolds.
\bpr\label{p:p}
    If $\pi:\mathbb CH^m_s\to B^n$ is a semi-Riemannian 
    submersion with complex, connected,
    totally geodesic fibres then $2m=n+2s$, 
    the  induced  metrics on fibres are 
    negative definite and the fibres are diffeomorphic to $\mathbb CP^s$.
\epr
\begin{proof}
    We denote by $J$ the natural almost complex structure 
    on $\mathbb CH^m_s$.
    By lemma \ref{lemtrei}, we have\\
a) $\hat R(U,V,U,V)=R(U,V,U,V)=-(g(U,U)g(V,V)-g(U,V)^2+3g(U,JV)^2)$.
   Hence the fibres have constant holomorphic curvature $-4$.\\
b) $g(A_XU,A_XU)=-(g(U,U)g(X,X)+3g(X,JU)^2)=-g(U,U)g(X,X)$, since the fibres are
   complex submanifolds. 
   We obtain $g(U,U)\leq 0$ for every vertical vector field $U$.\\
c) $R'(\pi_*X,\pi_*Y,\pi_*X,\pi_*Y)=R(X,Y,X,Y)+3g(A_XY,A_XY)=$\\
   $-(g(X,X)g(Y,Y)-g(X,Y)^2+3g(X,JY)^2)
   +3g(A_XY,A_XY)\leq 0$, since the induced metrics 
    on fibres are negative definite.
   By proposition \ref{simply}, it follows that the 
   fibres are simply connected. Since the fibres are 
   complete, simply connected, complex manifolds 
   with constant holomorphic curvature $-4$,
   we have that the fibres are isometric to $\mathbb CH^s_s$.
\end{proof} 
\bth\label{cgeo}
   If $\pi:\mathbb CH^m_s\to B$ is a semi-Riemannian submersion 
   with connected, complex, totally geodesic fibres  
   from a complex
   pseudo-hyperbolic space, then $\pi$
   is, up to equivalence, the canonical semi-Riemannian submersion
   given by example $4$
    $$\mathbb CH^{2k+1}_1\to\mathbb HH^k.$$
\eth
\begin{proof}
    Let $\theta:H^{2m+1}_{2s+1}\to\mathbb CH^m_s$
    be the canonical semi-Riemannian
    submersion with totally geodesic fibres given in 
    the definition \ref{d:3} (see also \cite{bar}
    or \cite{kwo}).
   We have $\tilde\pi=\pi\circ\theta:H^{2m+1}_{2s+1}\to B$ 
   is a semi-Riemannian 
   submersion with totally geodesic fibres, by theorem 2.5 in \cite{esco}.
   Since the dimension of fibres of $\tilde\pi$ is greater than or equal to 2,
   we get, by main theorem \ref{main}, the following
   possible situations:
\begin{itemize}
	\item[$i)$]  $m=2k+1$, $2s+1=3$ and $B$ is isometric to $\mathbb HH^k$ or
	\item[$ii)$] $m=7$, $2s+1=7$ and $B$ is isometric to $H^8(-4)$.
\end{itemize}
First, we shall prove that  any two semi-Riemannian submersions
   $\pi,\pi':\mathbb CH^{2k+1}_1\to\mathbb HH^k$ 
   with connected, complex, totally geodesic fibres 
   are equivalent.
 
By proof of proposition \ref{p:p}, we have
$g(A_XU,A_XU)=-g(U,U)g(X,X)$. 
Let $p,q\in\mathbb  CH^{2k+1}_1$
By proof of the main theorem, this implies 
$A^vA^w+A^wA^v=-2\tilde g(v,w)Id$. 
The extension of 
$\mathcal U:\mathcal V_p\to End(\mathcal H_p)$ constructed in proof of the main theorem,
to the Clifford algebra $Cl(\mathcal V_p, \tilde g_p)$
makes $\mathcal H_p$ a $Cl(\mathcal V_p, \tilde g_p)$-module
which splits in $k$ irreducible modules of dimension $4$.
By classification of irreducible representation for case $\dim\mathcal V_p=2$ and
$\tilde g_p$ positive definite, we have
any two such irreducible $Cl(\mathcal V_p, \tilde g_p)$-modules are equivalent.
Like in proof of the main theorem,
we may construct an isometry 
$L=L_1\oplus L_2:T_q\mathbb CH^{2k+1}_1\to T_p\mathbb CH^{2k+1}_1$, which
maps $\mathcal H'_q$ onto $\mathcal H_p$ and $A'$ onto $A$. This produces
an isometry $f:\mathbb CH^{2k+1}_1\to\mathbb CH^{2k+1}_1$
with $f(q)=p$ and $f_{*q}=L$ (see corollary 2.3.14 in \cite{wol}).
Again by theorem \ref{t:esc}, we get $\pi$ and $\pi'$ are equivalent.

For the case $ii)$ we shall obtain that there are no
   $\pi:\mathbb CH^7_3\to H^8(-4)$
   semi-Riemannian submersions with complex, connected, totally geodesic fibres.

\bpr\label{nocgeo}
    There are no $\pi:\mathbb CH^7_3\to H^8(-4)$ 
    semi-Riemannian submersions with connected, complex, totally
    geodesic fibres. 
\epr
\begin{proof}
%[Sketch of proof]
   The proof is based on Ranjan's argument (see proof of main theorem in \cite{ran}).
   Here, we show how to modify  Ranjan's argument to our different situation.

   Suppose there is $\pi:\mathbb CH^7_3\to H^8(-4)$
   a semi-Riemannian submersion with complex, connected, 
   totally geodesic fibres. By main theorem \ref{main},
   $\tilde\pi=\pi\circ\theta:H^{15}_7\to H^8(-4)$ is equivalent to the 
   canonical semi-Riemannian submersion
   $Spin(1,8)/Spin(7)\to Spin(1,8)/Spin(8)$ given by example 3.\\
   Let $\sigma:Spin(1,8)\to SO(8,8)$ be the spin representation
   of $Spin(1,8)$. $Spin(1,8)$ acts
   on $H^8(-4)$ via double covering map $Spin(1,8)\to SO(1,8)$
   and transitively on 
   $H^{15}_7\subset\mathbb R^{16}_8$.
   We denote by $Cl^0(\mathbb R^9_1)$ the even component of Clifford algebra 
   $Cl(\mathbb R^9_1)$.
   Notice that $Cl^0(\mathbb R^9_1)\cong M(16,\mathbb R)$,
   $Cl(\mathbb R^9_1)\cong M(16,\mathbb R)\oplus M(16,\mathbb R)$
    and the volume element $\omega$ in $Cl(\mathbb R_1^9)$ satisfies
    $\omega^2=1$ (see \cite{law}).

    For any $b\in H^8(-4)$, let $G_b$ be the isotropy group of $b$ in $Spin(1,8)$. 
    If we restrict $\sigma\left|_{G_b}\right.$ then  $\sigma\left|_{G_b}\right.$ 
    breaks $\mathbb R^{16}_8$ into two $\f 12$-spin representations. We will denote them
    by $\mathbb R^8_{\pm}$. Hence $\mathbb R^8_+\cap H^{15}_7=\tilde\pi^{-1}(b)$.
    Let $b^\perp=\{x\in\mathbb R^9_1\ |\ <x,b>=0\}$. 
    We have $Cl(b^\perp)\cap Spin(1,8)=G_b$, $\dim b^\perp=8$ and the following diagram
    is commutative
$$\begin{CD}
G_b @>>> Cl^0(b^\perp)\\
@VVV @VVV\\
Spin(1,8) @>>> Cl^0(\mathbb R^9_1),
\end{CD}
$$
where all arrows are standard inclusions.
Let $\{e_1,\dots, e_8\}$ be an orientated basis of $b^\perp$. Then $z'=e_1\dots e_8$
lies in the centre of $Cl^0(b^\perp)$ and $z'$ acts by $Id$ on $\mathbb R^8_+$ 
and $-Id$ on $\mathbb R^8_-$. We have $Cl(\sigma)(z')=\pm 1$ on $\mathbb R^8_\pm$.

Since $\mathbb R^8_+\cap H^{15}_7=\tilde\pi^{-1}(b)$, $\mathbb R^8_+$ is invariant under
$J$ and so is $\mathbb R^8_-$.
Here $J$ denotes the natural complex structure on $\mathbb R^{16}=\mathbb C^8$.
Hence $Cl(\sigma)(z')$ commutes with $J$.
Let $z\in Cl(\mathbb R^9_1)$ be the generator of the center of $Cl(\mathbb R^9_1)$.
We have either $z=e_1e_2\dots e_8b$ or $zb=-e_1e_2\dots e_8=-z'$.
Therefore $Cl(\sigma)(zb)$ commutes with $J$ for every $b\in H^8(-4)$
and hence for every $b\in\mathbb R^9$.

Consider the linear map $\alpha:\mathbb R^9\to M(16,\mathbb R)$ given by
$b\mapsto Cl(\sigma)(zb)$. It has the following properties:
\begin{itemize}
\item[$i)$]
   It factors through $M(8,\mathbb C)\subset M(16,\mathbb R)$
\item[$ii)$]
    $[Cl(\sigma)(zb)]^2=Cl(\sigma)((zb)^2)=Cl(\sigma)(-|b|^2)=-|b|^2Id.$
\end{itemize}
Hence $\alpha$ extends to a homomorphism
	$Cl(\alpha):Cl(\mathbb R^9_1)\to M(8,\mathbb C)$.
But $Cl(\mathbb R^9_1)\cong M(16,\mathbb R)\oplus M(16,\mathbb R)$ (see \cite{law}).
So the above homomorphism is impossible to exist.
We get the required contradiction.
\end{proof}
  This ends the proof of theorem \ref{cgeo}.\end{proof}
%%%%%%%%%%%%%%%%%%%%%%%%%%%%%%%%%%%%%%%%%%%%%%%%%%%%%%%%%%%%
\begin{acknow}
We would like to thank Dmitri Alekseevsky
for the remarks and comments  made to an earlier
version of this work.
The second author thanks to Stefano Marchiafava
for useful discussions on this topic with 
the occasion of his visit to Rome
in autumn of 1999.
\end{acknow}
%
% 
%\newpage
%\small
%\footnotesize
%\bibliographystyle{plain}

\end{document}